\documentclass[12pt,reqno]{amsart}

\usepackage{geometry}
\usepackage[T1]{fontenc}
\usepackage{lmodern,microtype}
\usepackage{amsmath,amssymb,amsthm,mathtools,mathrsfs}
\usepackage{fancyhdr,hyperref}
\hypersetup{hidelinks}

\newcommand{\sym}{\operatorname{Sym}}
\newcommand{\GL}{\mathrm{GL}}
\newcommand{\SL}{\mathrm{SL}}
\newcommand{\e}{\mathrm e}
\newcommand{\eps}{\varepsilon}
\newcommand{\R}{\mathbb R}
\newcommand{\Z}{\mathbb Z}

\newcommand{\GR}{\Gamma_{\mathbb R}}
\newcommand{\GC}{\Gamma_{\mathbb C}}
\newcommand{\intc}[1]{\frac{1}{2\pi i}\int_{(#1)}}
\newcommand{\dd}{\mathop{}\!\mathrm d}
\newtheorem{theorem}{Theorem}[section]
\newtheorem{proposition}[theorem]{Proposition}
\newtheorem{lemma}[theorem]{Lemma}
\newtheorem{corollary}[theorem]{Corollary}
\theoremstyle{remark}
\newtheorem{remark}{Remark}[section]
\numberwithin{equation}{section}

\title[The first moment of symmetric cube $L$-functions]
{The first moment of symmetric cube $L$-functions}
\author{Bingrong Huang}
\address{Data Science Institute \& School of Mathematics, Shandong University, Jinan, Shandong 250100, China} 
\curraddr{Department of Mathematics, Princeton University, Fine Hall, Washington Road, Princeton, NJ 08544, USA}
\email{brhuang@sdu.edu.cn}

\author{Anders S\"odergren}
\address{Department of Mathematical Sciences, Chalmers University of Technology and the University of Gothenburg, SE-412 96 Gothenburg, Sweden} 
\email{andesod@chalmers.se}

\date{\today}

\thanks{B.H. was supported by the Key R\&D Program of Shandong Province, China (No. 2026ZLYS01) and the Scientific Research Innovation Capability Support Project for Young Faculty (No. SRICSPYF-ZY2025158). A.S. was supported by the grant 2021-04605 from the Swedish
Research Council.}

\begin{document}

\begin{abstract}
We establish a power saving asymptotic formula for the first moment of symmetric cube central $L$-values over Hecke--Maass forms in short spectral windows.
As applications, we prove quantitative non-vanishing, an almost-all Lindel\"of bound, and almost-all decay of cubic periods.
We also prove the corresponding first moment asymptotic for central derivatives in the holomorphic weight aspect.

Our proof begins with an application of the Kuznetsov formula, followed by a first Poisson summation and stationary phase analysis. The key observation is that a periodic cubic factor can be separated from the resulting oscillation and incorporated into the arithmetic sum, allowing a second application of Poisson summation to yield further cancellation.
\end{abstract}
\maketitle

\section{Introduction}

Moments of central values of $L$-functions provide information on the distribution and typical size of $L$-values in families, with applications to problems such as non-vanishing and subconvexity. Spectral families of automorphic forms provide a natural setting for these questions.

Of particular interest are the symmetric power $L$-functions attached to $\GL(2)$ forms. The symmetric square lift is automorphic on $\GL(3)$ by Gelbart--Jacquet \cite{GJ}, and its central values have been studied extensively; see, for example, \cite{BF18,BF21,BF26a,BF26b,Blomer,Khan07,Khan,KY,KSeng,Lam,Lau,Liu,LS95,TX}. 
Among these works, Luo--Sarnak \cite{LS95} established an early bound on the second moment in the Maass form spectral aspect. In the holomorphic weight aspect, Lau \cite{Lau} obtained an asymptotic formula for the first moment, while Khan \cite{Khan} established an asymptotic formula for the twisted second moment after averaging over the weight.
By comparison, far fewer results are known for central values of symmetric cube $L$-functions.
The symmetric cube lift is automorphic on $\GL(4)$ by Kim--Shahidi \cite{KS02}. Another important result in this setting is due to Hoffstein--Jung--Lee \cite{HJL}, who proved that infinitely many Hecke--Maass cusp forms over $\mathbb Q(\sqrt{-3})$ have non-vanishing central values of their symmetric cube $L$-functions.

In this paper, we study the first moment of symmetric cube central $L$-values over Hecke--Maass forms for $\SL_2(\Z)$ in a short spectral window of length $M$. We prove an asymptotic formula with a power saving error term uniformly for $T^{3/4+\delta}\le M\le T^{1-\delta}$, where $0<\delta<1/8$ is fixed. We also derive quantitative non-vanishing, an almost-all Lindel\"of bound, and an almost-all bound for cubic periods.

\subsection{Main result}

Let $\mathcal B$ be an orthonormal Hecke basis of real-valued weight zero Hecke--Maass cusp forms for $\SL_2(\Z)$, including both even and odd forms. For $\phi\in\mathcal B$, let $t_\phi\ge0$ denote its spectral parameter, so that the Laplace eigenvalue is $1/4+t_\phi^2$. Moreover, normalize the Hecke eigenvalues by $\lambda_\phi(1)=1$. If $\alpha_\phi(p)$ and $\beta_\phi(p)$ are the Satake parameters at $p$, with $\alpha_\phi(p)\beta_\phi(p)=1$, then for $m=2,3$ and $\Re(s)>1$, we write
\[
  L(s,\sym^m\phi)
  =\prod_p\prod_{j=0}^{m}
  \left(1-\alpha_\phi(p)^{m-j}\beta_\phi(p)^j p^{-s}\right)^{-1}.
\]
Note that with this normalization, the functional equation is centered at $s=1/2$.

Throughout, $h$ is an even entire function, real-valued on $\R$, such that for every $A,B\ge0$ and every integer $j\ge0$,
\[
  h^{(j)}(x+iy)\ll_{A,B,j,h}(1+|x|)^{-A},
  \qquad |y|\le B.
\]
In particular, the Gaussian $h(x)=e^{-x^2}$ is admissible. 
%Fix $0<\delta<1/8$ and assume that
%\[
%  T^{3/4+\delta}\le M\le T^{1-\delta}.
%\]
%The letter $\eps$ denotes an arbitrarily small positive number whose value may vary from one occurrence to another.

\begin{theorem}\label{thm:main}
Let $0<\delta<1/8$. As $T\to\infty$, uniformly for
$T^{3/4+\delta}\le M\le T^{1-\delta}$, we have
\begin{equation}\label{eq:main}
  \sum_{\phi\in\mathcal B}
   \frac{h((t_\phi-T)/M)}{L(1,\sym^2\phi)}
   L\!\left(\tfrac12,\sym^3\phi\right)
   =\frac13 TM\int_{\R} h(x)\,\dd x
   +O_{h,\delta,\eps}\!\left(\frac{T^{5/2+\eps}}M\right),
\end{equation}
for any $\eps>0$. 
% In particular, the error term is power-saving throughout the stated range.
\end{theorem}

\begin{remark}
For comparison, when $t_\phi\asymp T$, the symmetric square $L$-function has analytic conductor $\asymp T^2$, whereas $L(s,\sym^3\phi)$ has analytic conductor $\asymp T^4$ in the sense of Iwaniec--Sarnak \cite{IS}. Thus the corresponding approximate functional equations have natural lengths about $T$ and $T^2$, respectively. The larger conductor in the symmetric cube case leads to a longer Dirichlet polynomial and is one source of difficulty in obtaining an asymptotic spectral first moment.
\end{remark}

\begin{remark}
The argument for Theorem~\ref{thm:main} also extends to first moments twisted by a Hecke eigenvalue $\lambda_\phi(m)$, with at most polynomial dependence on $m$ in the error term. Since Theorem~\ref{thm:main} has a power saving, the twisting parameter $m$ may grow as a sufficiently small fixed power of $T$. Such twisted first moment formulas can be used, for example, to remove the harmonic weights in Theorem~\ref{thm:main} and, more generally, as input for lower bounds for higher moments of $L(1/2,\sym^3\phi)$. We do not pursue these extensions here.
\end{remark}

The proof begins with an exact Dirichlet series identity expressing the symmetric cube coefficients in terms of ordinary Hecke eigenvalues:
\begin{align}\label{eq:coefficient identity}
  L(s,\sym^3\phi)=\sum_{n\ge1}\frac{A_\phi(n)}{n^s},
  \qquad
  A_\phi(n)=\sum_{n_1n_2^2n_3^3n_4^4=n}\mu^2(n_3)
  \lambda_\phi(n_1^3n_2^2n_3^3).
\end{align}
This identity allows us to apply the Kuznetsov formula directly to the first moment. The diagonal is then given by $n_1=n_2=n_3=1$, which corresponds to the fourth power terms in the original Dirichlet series.

The main difficulty lies in estimating the resulting Kloosterman terms. Our approach uses two successive Poisson summations. After the first Poisson summation, stationary phase produces the two oscillatory factors
\[
  \e\!\left(-\frac{m^3}{27n_2^2n_3^3c}\right)
  \e\!\left(-\frac{27cn_2^2n_3^3T^2}{4\pi^2m^3}\right).
\]
The first factor, although arising from the analytic transform, is periodic in $m$ modulo $27n_2^2n_3^3c$ and can therefore be incorporated into the arithmetic sum, while the second remains in the analytic weight. This analytic-to-arithmetic transition leaves a more slowly oscillating analytic weight, allowing the second Poisson summation to yield further cancellation. Related interactions between analytic oscillation and arithmetic characters can be found in Conrey--Iwaniec \cite[\S 8]{CI} and Li \cite[\S 4]{Li}.

\subsection{Applications}

By Lapid--Rallis \cite[Theorem~1]{LR}, the central values in Theorem~\ref{thm:main} are nonnegative; see also Lemma~\ref{lem:fe}. We record three consequences of Theorem~\ref{thm:main}.
Throughout, let $0<\delta<1/8$ and
$T^{3/4+\delta}\le M\le T^{1-\delta}$.

\begin{corollary}[Quantitative non-vanishing]\label{cor:nonvan}
We have
\begin{equation}\label{eq:quantnonvan}
  \#\{\phi\in\mathcal B:|t_\phi-T|\le M,\ L(\tfrac12,\sym^3\phi)\ne0\}
  \gg_{\delta,\eps}M^{10/9}T^{-\eps}.
\end{equation}
\end{corollary}

The exponent $10/9$ comes from the tenth-moment estimate of Huang--Li \cite[Theorem~1.4]{HL}. By the Eisenstein series criterion of Lapid--Rallis \cite[p.~892]{LR}, every form counted in Corollary~\ref{cor:nonvan} also gives a simple pole at $s=1/2$ of the Siegel Eisenstein series induced from $\sym^3\pi_\phi$, whose residue is square integrable.

\begin{corollary}[Almost-all Lindel\"of]\label{cor:lindelof}
For $Y\ge1$,
\[
  \#\{\phi\in\mathcal B:|t_\phi-T|\le M,\ L(\tfrac12,\sym^3\phi)>Y\}
  \ll_{\delta,\eps}\frac{TM T^\eps}{Y}.
\]
Thus, for each fixed $a>0$, all but $O_{\delta,a}(TM T^{-a/2})$ forms in the spectral window $|t_\phi-T|\le M$ satisfy
\[
  L(\tfrac12,\sym^3\phi)\le T^a.
\]
\end{corollary}

The first moment asymptotic in Theorem \ref{thm:main} also gives a bound for cubic periods of Maass forms. Write
\[
  \langle1,\phi^3\rangle
  =\int_{\SL_2(\Z)\backslash\mathbb H}
     \phi(z)^3\,\frac{dx\,dy}{y^2}.
\]

\begin{corollary}[Almost-all cubic periods]\label{cor:cubic}
We have
\begin{equation}\label{eq:cubic-first}
  \sum_{T\le t_\phi\le T+M}|\langle1,\phi^3\rangle|
  \ll_{\delta,\eps}MT^\eps.
\end{equation}
Consequently, for each fixed $a>0$, all but
$O_{\delta,a}(TM\,T^{-a/2})$ forms in this interval satisfy
\[
  |\langle1,\phi^3\rangle|\le T^{-1+a}.
\]
\end{corollary}

In the present range of $M$, this improves the almost-all decay bound following from the variance estimate of Huang--Li \cite[Theorem~1.1]{HL}, which, however, applies to much shorter spectral intervals.

\subsection{A holomorphic analogue}

In the holomorphic setting, the functional equation of the symmetric cube $L$-function has sign $-1$, so the natural central quantity is the derivative. Let $H_k$ be a Hecke basis of $S_k(\SL_2(\Z))$. 

\begin{theorem}\label{thm:holo}
Let $\Psi\in C_c^\infty(\R)$ be real-valued and
$0<\delta<1/8$. As $K\to\infty$, uniformly for
$K^{3/4+\delta}\le M\le K^{1-\delta}$, we have
\begin{equation}\label{eq:holomain}
\begin{split}
 &\sum_{\substack{k\ge2\\k\equiv0\bmod 2}}\Psi\!\left(\frac{k-K}{M}\right)
       \sum_{f\in H_k}\frac{L'(1/2,\sym^3 f)}{L(1,\sym^2 f)}\\
 &\quad=\sum_{\substack{k\ge2\\k\equiv0\bmod 2}}\Psi\!\left(\frac{k-K}{M}\right)
       \frac{k-1}{\pi^2}
       \Bigg\{\zeta(2)\bigg[\psi\!\left(\frac{k}{2}\right)
          +\psi\!\left(\frac{3k-2}{2}\right)-2\log(2\pi)\bigg]
          +4\zeta'(2)\Bigg\}\\
 &\qquad\qquad+O_{\Psi,\delta,\eps}\!\left(\frac{K^{5/2+\eps}}M\right),
\end{split}
\end{equation}
where $\psi=\Gamma'/\Gamma$. 
The right-hand side of \eqref{eq:holomain} is asymptotic to $\tfrac16(KM\log K)\int_{\R} \Psi(x)\,\dd x$.
\end{theorem}

The exact main sum depending on $k$ is retained to avoid hiding a weight window Taylor error term in the stated remainder. If $\Psi$ is supported in $[-1,1]$ and $\int_{\R} \Psi(x)\,\dd x>0$, Theorem~\ref{thm:holo} implies that $L'(1/2,\sym^3 f)\ne0$ for some even weight $|k-K|\le M$, for every sufficiently large $K$. 
%The proof is given in the final section.

\subsection{Sketch of the proof of Theorem~\ref{thm:main}}
\label{subsec:sketch}

To illustrate the main mechanism without obscuring it by the outer divisor variables, we consider the model case $n_2=n_3=n_4=1$, write $n=n_1\asymp N$, and suppress harmless smooth weights. After the approximate functional equation, the coefficient identity \eqref{eq:coefficient identity}, and the Kuznetsov formula, a typical Kloosterman block has the shape
\[
  \mathscr O(N,C)=TM\sum_{c\asymp C}\frac1{\sqrt c}
  \sum_{n\asymp N}\frac{S(1,n^3;c)}{n^{5/4}}W(n,c)
  \e\!\left(
    \frac{2n^{3/2}}c
    -\frac{cT^2}{4\pi^2n^{3/2}}
  \right),
\]
where the Bessel transform restricts
\[
  C\ll T^\eps\frac{N^{3/2}}{TM}.
\]

We first apply Poisson summation to the $n$-sum modulo $c$. This introduces the complete sum
\[
  \mathfrak C_c(m)
  =\sum_{x\bmod c}S(1,x^3;c)\e_c(mx)
\]
and an oscillatory integral in the dual variable $m$. Stationary phase restricts
\[
  m\asymp N^{1/2}
\]
and, up to a smooth amplitude, splits the oscillation into
\[
  \e\!\left(-\frac{m^3}{27c}\right)
  \e\!\left(-\frac{27cT^2}{4\pi^2m^3}\right).
\]
At this stage the block has, up to a negligible error term, the schematic form
\[
  \frac{TM}{N}
  \sum_{c\asymp C}\frac1c
  \sum_{m\asymp N^{1/2}}
  \mathfrak C_c(m)W_1(m,c)
  \e\!\left(-\frac{m^3}{27c}
             -\frac{27cT^2}{4\pi^2m^3}\right).
\]
However, this is not yet sufficient. Indeed, even granting square root cancellation
\[
  \mathfrak C_c(m)\ll_\eps c^{1+\eps},
\]
estimating the $m$- and $c$-sums trivially gives only
\[
  \frac{TM}{N}\,C N^{1/2}T^\eps
  \ll N^{1+\eps},
\]
where we used $C\ll T^\eps N^{3/2}/(TM)$. At the largest scale $N\asymp T^2$, this is $T^{2+\eps}$, which is larger than the main term of size $TM$ by essentially a factor $T/M$. Thus, the first Poisson summation, stationary phase, and square root cancellation in the resulting complete sums do not by themselves give a power saving; an additional saving of essentially $T/M$ is needed.

We now incorporate the first factor, which is periodic in $m$ modulo $27c$, into the arithmetic sum and retain the second factor in the analytic weight. The effect of this separation is clearest at the largest relevant scale. When
\[
  N\asymp T^2,
  \qquad
  C\asymp \frac{T^2}{M},
  \qquad
  m\asymp T,
\]
the two phases have respective sizes
\[
  \frac{m^3}{27c}\asymp TM,
  \qquad
  \frac{27cT^2}{4\pi^2m^3}\asymp\frac{T}{M}.
\]
Thus the large analytic oscillation of size $TM$ has become arithmetic, while only an oscillation of size $T/M$ remains in the analytic weight.

We next apply Poisson summation to the $m$-sum modulo $27c$. The resulting complete sums have the form
\[
  \mathfrak S(\ell;c)
  =\sum_{m\bmod 27c}\mathfrak C_c(m)
   \e_{27c}(-m^3+\ell m).
\]
At the largest scale, the second dual variable is restricted to
\[
  |\ell|\asymp \frac{T^2}{M^2},
\]
so the resulting dual sum is much shorter than the original sum. For a generic prime modulus $c=p$, $p\ge5$, the complete sum exhibits square root cancellation. For $x\ne0\pmod p$, we use the identity
\[
  S(1,x^3;p)=S(x,x^2;p).
\]
This makes the dependence on $x$ quadratic, allowing two quadratic Gauss summations and reducing the local factor to a one-variable mixed character sum. Consequently, the local complete sum is $O(p^{3/2})$, a square root saving over the trivial bound. The losses from the remaining prime powers are recorded by a powerful part of the modulus and are summable on average. Finally, together with the analysis of the second integral transform, these estimates give, in the model case,
\[
  \mathscr O(N,C)
  \ll_\eps T^\eps\left(MN^{1/4}+\frac{N^{5/4}}M\right).
\]
At $N\asymp T^2$, this becomes
\[
  \mathscr O(N,C)
  \ll_\eps T^\eps\left(MT^{1/2}+\frac{T^{5/2}}M\right)
  \ll_\eps \frac{T^{5/2+\eps}}M,
\]
which is a power saving over $TM$ throughout the range of Theorem~\ref{thm:main}.

%The holomorphic case follows the same architecture, with the Petersson formula replacing the Kuznetsov formula; the details are given in Section~\ref{sec:holo}.

\subsection*{Organization and notation}

Section~\ref{sec:prelim} collects the analytic and automorphic preliminaries, including the approximate functional equation, Poisson summation, and the Kuznetsov formula. Section~\ref{sec:proof} proves Theorem~\ref{thm:main}. Section~\ref{sec:applications} proves the three corollaries, and Section~\ref{sec:holo} treats the holomorphic case.

We write $\e(x)=e^{2\pi i x}$ and $\e_q(x)=\e(x/q)$. The letter $\eps>0$ denotes an arbitrarily small positive constant, whose value may vary from one occurrence to another, whereas $0<\delta<1/8$ is fixed. We use $A\ll B$ and $A=O(B)$ interchangeably, with subscripts indicating the permitted dependence of the implied constant. The implied constant is allowed to depend on the fixed test functions and on $\delta$ without explicit mention. For a prime $p$, $v_p(n)$ denotes the exponent of $p$ in $n$. If $(a,q)=1$, then $\bar a$ denotes the multiplicative inverse of $a$ modulo $q$.

\section{Preliminaries}\label{sec:prelim}

\subsection{Functional equation and coefficient bounds}
For $\phi\in\mathcal B$, write
\[
  \phi(z)=\sqrt y\sum_{n\ne0}\rho_\phi(n)K_{it_\phi}(2\pi|n|y)\e(nx),
  \qquad \rho_\phi(n)=\rho_\phi(1)\lambda_\phi(n)\quad(n>0),
\]
where $\lambda_\phi(1)=1$. Let $\kappa_\phi=0$ if $\phi$ is even and
$\kappa_\phi=1$ if $\phi$ is odd, so that
$\rho_\phi(-n)=(-1)^{\kappa_\phi}\rho_\phi(n)$ for $n>0$. Set
\begin{equation}\label{eq:gamma}
  \GR(s)=\pi^{-s/2}\Gamma(s/2),\qquad
  \gamma_\kappa(s,t)=\prod_{\nu=\pm1,\pm3}\GR(s+\kappa+i\nu t).
\end{equation}

\begin{lemma}\label{lem:fe}
The function $\Lambda(s,\sym^3\phi)=\gamma_{\kappa_\phi}(s,t_\phi)L(s,\sym^3\phi)$ is entire and satisfies
\[
  \Lambda(s,\sym^3\phi)=\Lambda(1-s,\sym^3\phi).
\]
%Its finite conductor is $1$, its root number is $+1$, and its central analytic conductor is $\asymp(1+t_\phi)^4$. 
Moreover, $L(1/2,\sym^3\phi)\ge0$.
\end{lemma}
\begin{proof}
Let $\pi_\phi$ be the automorphic representation generated by $\phi$.
The symmetric cube lift $\Pi=\sym^3\pi_\phi$ is cuspidal by
Kim--Shahidi \cite[Theorem~6.1]{KS02}. 
%Indeed, the exceptions in that
%criterion require a nontrivial quadratic self-twist of $\pi_\phi$ or a
%nontrivial cubic self-twist of $\sym^2\pi_\phi$. Since both representations
%are unramified at every finite prime, such a twisting character would
%also be unramified there. No nontrivial finite-order Hecke character of
%$\mathbb Q$ has this property. Thus the completed standard $L$-function
%of $\Pi$ is entire.
%
By \cite[Theorem~4.2(a),(b)]{KS00}, the local factors of $\Pi$ are the
symmetric cubes of those of $\pi_\phi$. At infinity, the four parameters
are $\kappa_\phi+i\nu t_\phi$, $\nu=\pm1,\pm3$, giving
\eqref{eq:gamma} and local root number $i^{4\kappa_\phi}=1$.
Cuspidality gives entireness, while self-duality and unramifiedness at
every finite prime give the stated functional equation.
%At every finite prime the lift is unramified, so its conductor exponent
%is zero and its local root number is $1$. The trivial central character
%of $\pi_\phi$ makes $\Pi$ self-dual, so these local calculations give the
%stated functional equation. Its central analytic conductor is
%comparable to
%\[
%  \prod_{\nu=\pm1,\pm3}
%  \bigl(1+|\kappa_\phi+i\nu t_\phi|\bigr)
%  \asymp(1+t_\phi)^4.
%\]
Finally, $\Pi$ is of symplectic type by \cite[Section~9]{KS02}, so
Lapid--Rallis \cite[Theorem~1]{LR} gives
$L(1/2,\sym^3\phi)\ge0$.
\end{proof}

For $\Re s$ sufficiently large, define the Dirichlet coefficients
$A_\phi(n)$ by
\[
  L(s,\sym^3\phi)=\sum_{n\ge1}\frac{A_\phi(n)}{n^s}.
\]
We shall use the following estimates, which are uniform in $\phi$.

\begin{lemma}\label{lem:coeffmean}
Let $T\ge2$ and $t_\phi\asymp T$. For $x\ge2$ and any fixed $B>2$,
\begin{align}
  \sum_{n\le x}\frac{|A_\phi(n)|}{n}&\ll_\eps(Tx)^\eps,\label{eq:coeffmean}\\
  \sum_{n\ge1}\frac{|A_\phi(n)|}{\sqrt n}(1+n/T^2)^{-B}
     &\ll_{B,\eps}T^{1+\eps}.\label{eq:coeffweighted}
\end{align}
\end{lemma}
\begin{proof}
Molteni's theorem \cite[Theorem~4]{Molteni} applies to the cuspidal
representation $\sym^3\pi_\phi$ on $\GL_4$; its hypotheses for cuspidal
automorphic $L$-functions are verified in the discussion following that
theorem. In the normalization \eqref{eq:gamma}, its main parameter is
$\asymp T^4$. Hence, for every $\eps>0$,
\[
  \sum_{n\le x}\frac{|A_\phi(n)|}{n}\ll_\eps (T^4x)^\eps,
\]
which gives \eqref{eq:coeffmean} after relabeling $\eps$.
To prove \eqref{eq:coeffweighted}, it suffices to take $0<\eps<1$.
Applying \eqref{eq:coeffmean} with $\eps/3$ on dyadic intervals gives
\[
  \begin{aligned}
  &\sum_{n\ge1}\frac{|A_\phi(n)|}{\sqrt n}(1+n/T^2)^{-B}
  \ll_\eps T^{\eps/3}
  \sum_{\substack{Z=2^j\\j\ge0}}
  Z^{1/2+\eps/3}(1+Z/T^2)^{-B}
  \ll_{B,\eps}T^{1+\eps},
  \end{aligned}
\]
where the last bound follows by splitting the dyadic sum at $Z=T^2$.
\end{proof}

\subsection{An exact Hecke expansion}\label{subsec:coeff}

\begin{lemma}\label{lem:coeffidentity}
For every $n\ge1$,
\begin{equation}\label{eq:single}
  A_\phi(n)=\sum_{n_1n_2^2n_3^3n_4^4=n}\mu^2(n_3)\lambda_\phi(n_1^3n_2^2n_3^3).
\end{equation}
\end{lemma}
\begin{proof}
Both sides of \eqref{eq:single} are multiplicative, so it suffices to
work at a prime $p$. The right-hand side has local generating series
\begin{equation}\label{eq:local}
  \begin{split}
  &\sum_{j\ge0}x^j
  \sum_{\substack{a,b,c,d\ge0\\a+2b+3c+4d=j}}
  \mu^2(p^c)\lambda_\phi(p^{3a+2b+3c})\\
  &\quad=\frac1{1-x^4}\sum_{a,b\ge0}x^{a+2b}
  \bigl(\lambda_\phi(p^{3a+2b})
       +x^3\lambda_\phi(p^{3a+2b+3})\bigr).
  \end{split}
\end{equation}
Let $z,z^{-1}$ be the Satake parameters at $p$. For $z\ne\pm1$, we have
\[
  \lambda_\phi(p^j)=\frac{z^{j+1}-z^{-j-1}}{z-z^{-1}},
  \qquad j\ge0.
\]
Summing the geometric series in \eqref{eq:local} therefore gives
\[
  \frac1{(z-z^{-1})(1-x^4)}
  \left\{
  \frac{z+x^3z^4}{(1-xz^3)(1-x^2z^2)}
  -\frac{z^{-1}+x^3z^{-4}}{(1-xz^{-3})(1-x^2z^{-2})}
  \right\}.
\]
Combining the fractions, we factor the numerator as
\[
  \begin{aligned}
  (z+x^3z^4)&(1-xz^{-3})(1-x^2z^{-2})-(z^{-1}+x^3z^{-4})(1-xz^3)(1-x^2z^2)\\
  &\quad=(z-z^{-1})(1+x^2-x^4-x^6)
        +(z^2-z^{-2})(x-x^5)\\
  &\quad=(z-z^{-1})(1-x^4)(1+xz)(1+xz^{-1}).
  \end{aligned}
\]
After cancellation, the series in \eqref{eq:local} becomes
\[
  \begin{aligned}
  &\frac{(1+xz)(1+xz^{-1})}
  {(1-xz^3)(1-x^2z^2)(1-xz^{-3})(1-x^2z^{-2})}
  =\prod_{\nu=\pm1,\pm3}(1-xz^\nu)^{-1}
  =\sum_{j\ge0}A_\phi(p^j)x^j,
  \end{aligned}
\]
where the last equality follows from the Euler product for
$L(s,\sym^3\phi)$. The local generating series and the Euler factor have
Laurent polynomial coefficients in $z$, so their equality as formal power
series in $x$ extends to $z=\pm1$. Comparing coefficients and using
multiplicativity proves \eqref{eq:single}.
\end{proof}

\subsection{Approximate functional equation}\label{subsec:afe}
For $n\ge1$ and $|\Im t|<1$, define
\begin{equation}\label{eq:weight}
  V(n,t)=\intc{3}\frac{e^{w^2}}w
    \left(\frac{3(1+t^2)}{4\pi^2n}\right)^w\,dw,
\end{equation}
where the logarithm is chosen real for real $t$.

\begin{lemma}\label{lem:afe}
Uniformly for $t_\phi\asymp T$,
\begin{equation}\label{eq:afe}
  L\!\left(\tfrac12,\sym^3\phi\right)
  =2\sum_{n\ge1}\frac{A_\phi(n)}{\sqrt n}V(n,t_\phi)
   +O_\eps(T^{-1+\eps}).
\end{equation}
For fixed $j,b,B\ge0$,
\begin{equation}\label{eq:weightbounds}
  (n\partial_n)^j\partial_t^bV(n,t)
  \ll_{j,b,B,\eps}T^{-b+\eps}(1+n/T^2)^{-B},
  \qquad t\asymp T.
\end{equation}
\end{lemma}
\begin{proof}
The approximate functional equation \cite[Theorem~5.3]{IK} gives
\[
  L\!\left(\tfrac12,\sym^3\phi\right)
  =2\sum_{n\ge1}\frac{A_\phi(n)}{\sqrt n}
       V_{\kappa_\phi}^{*}(n,t_\phi),
\]
where
\[
  V_\kappa^{*}(n,t)=\intc{3}\frac{e^{w^2}}w
    \frac{\gamma_\kappa(1/2+w,t)}{\gamma_\kappa(1/2,t)}n^{-w}\,dw.
\]
Stirling's formula for the pairs with imaginary parts $\pm t$ and $\pm3t$
gives, uniformly for $\kappa=0,1$ on fixed vertical lines with
$|\Im w|\le T^\eps$,
\[
  \frac{\gamma_\kappa(1/2+w,t)}{\gamma_\kappa(1/2,t)}
  =\left(\frac{3(1+t^2)}{4\pi^2}\right)^w
   \left\{1+O_\eps\!\left(T^{-2+\eps}(1+|w|)^C\right)\right\}.
\]
The Gaussian makes the remaining tails negligible. The gamma quotient
and its leading term both equal $1$ at $w=0$, so their difference divided
by $w$ is holomorphic there. Shifting the contour for
$V_\kappa^{*}(n,t)-V(n,t)$ to $\Re w=0$ when $n\le T^2$, and to a
sufficiently large fixed positive line otherwise, yields
\[
  V_\kappa^{*}(n,t)-V(n,t)
  \ll_{B,\eps}T^{-2+\eps}(1+n/T^2)^{-B}.
\]
By Lemma~\ref{lem:coeffmean}, replacing $V_{\kappa_\phi}^{*}$ by $V$
therefore costs $O_\eps(T^{-1+\eps})$, proving \eqref{eq:afe}.
Differentiation and contour shifting in \eqref{eq:weight} give
\eqref{eq:weightbounds}.
\end{proof}

In particular, \eqref{eq:afe} and Lemma~\ref{lem:coeffmean} give
\[
  L\!\left(\tfrac12,\sym^3\phi\right)\ll_\eps T^{1+\eps}.
\]

\subsection{Poisson summation}\label{subsec:poisson-prelim}
We use the Fourier transform convention
\[
  \widehat F(\xi)=\int_{\R}F(x)\e(-\xi x)\,\dd x.
\]
\begin{lemma}[Poisson summation {\cite[(4.24)]{IK}}]\label{lem:poisson}
Let $q\ge1$ be an integer, $B$ a $q$-periodic function on $\Z$, and $F$
a Schwartz function on $\R$. Then
\begin{equation}\label{eq:poisson}
  \sum_n B(n)F(n)=\frac1q\sum_{m\in\Z}
    \left(\sum_{b\bmod q}B(b)\e_q(mb)\right)\widehat F(m/q).
\end{equation}
\end{lemma}

\subsection{The Kuznetsov formula}\label{subsec:kuznetsov}
For positive integers $c$, write
\[
  S(m,n;c)=\sideset{}{^{*}}\sum_{d\bmod c}\e_c(md+n\bar d).
\]
We shall use the classical Weil bound \cite[Corollary~11.12]{IK}
\begin{equation}\label{eq:Weil}
  |S(m,n;c)|\le d(c)c^{1/2}(m,n,c)^{1/2}.
\end{equation}
Set
\[
  \tau_{it}(n)=\sum_{ab=n}(a/b)^{it},
  \qquad
  \alpha_\phi=\frac{|\rho_\phi(1)|^2}{\cosh(\pi t_\phi)}
     =\frac{2}{L(1,\sym^2\phi)}.
\]
The second identity follows from the Rankin--Selberg calculation with our Fourier normalization.

\begin{lemma}[Kuznetsov formula {\cite[Theorem~9.3]{IwaniecSpectral}}]\label{lem:kuznetsov}
Let $G$ be even, holomorphic in a fixed strip $|\Im t|\le 1/2+\sigma$
for some $\sigma>0$, and rapidly decreasing there. For integers $m,n\ge1$,
we have
\begin{equation}\label{eq:kuz}
\begin{split}
 &\sum_{\phi}\alpha_\phi\lambda_\phi(m)\lambda_\phi(n)G(t_\phi)
  +\frac1\pi\int_\R\frac{\tau_{it}(m)\tau_{it}(n)G(t)}{|\zeta(1+2it)|^2}\,\dd t\\
 &\quad=\frac{\delta_{m,n}}{\pi^2}\int_\R t\tanh(\pi t)G(t)\,\dd t
   +\sum_{c\ge1}\frac{S(m,n;c)}c
       \mathcal H_{G}\!\left(\frac{4\pi\sqrt{mn}}c\right),
\end{split}
\end{equation}
where
\begin{equation}\label{eq:besseldef}
  \mathcal H_{G}(x)=\frac{2i}\pi\int_\R
     \frac{tJ_{2it}(x)}{\cosh(\pi t)}G(t)\,\dd t.
\end{equation}
\end{lemma}

\section{Proof of Theorem~\ref{thm:main}}\label{sec:proof}

\subsection{Applying the approximate functional equation}\label{subsec:apply-afe}
For the trace formula we use the even spectral weight
\begin{equation}\label{eq:test}
  g(t)=h\!\left(\frac{t-T}{M}\right)
       +h\!\left(\frac{t+T}{M}\right).
\end{equation}
Write
\[
  \mathcal M_g=
  \sum_{\phi\in\mathcal B}
   \frac{g(t_\phi)}{L(1,\sym^2\phi)}
   L\!\left(\tfrac12,\sym^3\phi\right).
\]
By Lemma~\ref{lem:afe}, the standard weighted Weyl law
\[
  \sum_{\phi}\frac{|g(t_\phi)|}{L(1,\sym^2\phi)}
  \ll_{h,\eps} TM T^\eps,
\]
and the identity $\alpha_\phi=2/L(1,\sym^2\phi)$, we have
\begin{equation}\label{eq:afterafe}
  \mathcal M_g=\sum_{n_1,n_2,n_3,n_4\ge1}
  \frac{\mu^2(n_3)}{\sqrt{n_1n_2^2n_3^3n_4^4}}
  \sum_{\phi}\alpha_\phi\lambda_\phi(n_1^3n_2^2n_3^3)
  g(t_\phi)V(n_1n_2^2n_3^3n_4^4,t_\phi)
  +O_{h,\delta,\eps}(MT^\eps).
\end{equation}

\subsection{Applying the Kuznetsov formula}\label{subsec:apply-kuz}
The function $g$ is even and entire. Moreover, for every $A,B\ge0$ and every integer $j\ge0$,
\[
  g^{(j)}(t)\ll_{A,B,j,h}M^{-j}
  \left\{
    \left(1+\frac{|\Re t-T|}{M}\right)^{-A}
    +\left(1+\frac{|\Re t+T|}{M}\right)^{-A}
  \right\},
  \qquad |\Im t|\le B.
\]
By \eqref{eq:weight}, $V(n,t)$ is holomorphic with polynomial growth in the fixed strip $|\Im t|\le3/4$ for polynomially bounded $n$. Thus $g(t)V(n,t)$ is an admissible test function in \eqref{eq:kuz}.

Applying \eqref{eq:kuz} with $m=1$ and second Hecke index $n_1^3n_2^2n_3^3$ to \eqref{eq:afterafe}, we obtain
\begin{equation}\label{eq:decomposition}
  \mathcal M_g=\mathcal D+\mathcal O-\mathcal E+O(MT^\eps),
\end{equation}
where
\begin{align*}
  \mathcal D
  &=\frac1{\pi^2}\sum_{n_1,n_2,n_3,n_4\ge1}
     \frac{\mu^2(n_3)\delta_{n_1^3n_2^2n_3^3,1}}
          {\sqrt{n_1n_2^2n_3^3n_4^4}}
     \int_\R t\tanh(\pi t)g(t)V(n_1n_2^2n_3^3n_4^4,t)\,\dd t\\
  &=\frac1{\pi^2}\int_\R t\tanh(\pi t)g(t)
     \sum_{n_4\ge1}\frac{V(n_4^4,t)}{n_4^2}\,\dd t,
\end{align*}
\[
  \mathcal O=\sum_{n_1,n_2,n_3,n_4\ge1}
  \frac{\mu^2(n_3)}{\sqrt{n_1n_2^2n_3^3n_4^4}}
  \sum_{c\ge1}\frac{S(1,n_1^3n_2^2n_3^3;c)}c\,
  \mathcal H_{\,g(\cdot)V(n_1n_2^2n_3^3n_4^4,\cdot)}
  \!\left(\frac{4\pi n_1^{3/2}n_2n_3^{3/2}}c\right),
\]
and
\[
  \mathcal E=\frac1\pi\sum_{n_1,n_2,n_3,n_4\ge1}
  \frac{\mu^2(n_3)}{\sqrt{n_1n_2^2n_3^3n_4^4}}
  \int_\R\frac{\tau_{it}(n_1^3n_2^2n_3^3)}{|\zeta(1+2it)|^2}
  g(t)V(n_1n_2^2n_3^3n_4^4,t)\,\dd t.
\]
Thus $\mathcal D$, $\mathcal O$, and $\mathcal E$ are the diagonal, Kloosterman, and continuous spectrum terms, respectively.

\begin{lemma}\label{lem:diag-eis}
We have
\begin{align}
  \mathcal D
  &=\frac13 TM\int_\R h(x)\,\dd x
    +O_{h,\eps}(MT^{1/2+\eps}),\label{eq:diag}
\end{align}
and 
\begin{align}
  |\mathcal E|&\ll_{h,\eps}T^{1+\eps}.\label{eq:eis}
\end{align}
\end{lemma}
\begin{proof}
For the diagonal term,
\[
  \sum_{n_4\ge1}\frac{V(n_4^4,t)}{n_4^2}
  =\intc{3}\frac{e^{w^2}}w
    \left(\frac{3(1+t^2)}{4\pi^2}\right)^w
    \zeta(2+4w)\,dw.
\]
Shift the line of integration to $\Re w=-1/4+\sigma$, where $0<\sigma<1/4$. The only pole crossed is at $w=0$, with residue $\zeta(2)$, and hence, uniformly for $|t|\asymp T$,
\[
  \sum_{n_4\ge1}\frac{V(n_4^4,t)}{n_4^2}
  =\zeta(2)+O_\eps(T^{-1/2+\eps}).
\]
Therefore
\[
  \mathcal D
  =\frac{\zeta(2)}{\pi^2}
    \int_\R t\tanh(\pi t)g(t)\,\dd t
   +O_{h,\eps}(MT^{1/2+\eps}).
\]
Since $h$ is even,
\[
  \int_\R t\tanh(\pi t)g(t)\,\dd t
  =2TM\int_\R h(x)\,\dd x+O_{A,h,\delta}(T^{-A})
\]
for every $A>0$. Since $2\zeta(2)/\pi^2=1/3$, this proves \eqref{eq:diag}.

For the continuous spectrum, define $A_{it}(n)$ by
\[
  \prod_{\nu=\pm1,\pm3}\zeta(s+i\nu t)
  =\sum_{n\ge1}\frac{A_{it}(n)}{n^s},
  \qquad \Re(s)>1.
\]
The local identity in Lemma~\ref{lem:coeffidentity}, with $\lambda_\phi$ replaced by $\tau_{it}$, gives
\[
  \mathcal E
  =\frac1\pi\int_\R
    \frac{g(t)}{|\zeta(1+2it)|^2}
    \sum_{n\ge1}\frac{A_{it}(n)}{\sqrt n}V(n,t)\,\dd t.
\]
Repeating the contour argument in Lemma~\ref{lem:afe}, uniformly for $|t|\asymp T$, gives
\[
  2\sum_{n\ge1}\frac{A_{it}(n)}{\sqrt n}V(n,t)
  =\left|\zeta\!\left(\tfrac12+it\right)
          \zeta\!\left(\tfrac12+3it\right)\right|^2
   +O_\eps(T^{-1+\eps}).
\]
Using $|\zeta(1+2it)^{-1}|\ll\log(2+|t|)$ from \cite[\S1.43, Theorem~9*]{Titchmarsh}, together with the rapid decay of $g$, we obtain
\[
  |\mathcal E|
  \ll (\log T)^2
  \int_{T/2\le |t|\le 2T}
    \left|\zeta\!\left(\tfrac12+it\right)\right|^2
    \left|\zeta\!\left(\tfrac12+3it\right)\right|^2\,\dd t
  +O_{h,\eps}(MT^{-1+\eps}).
\]
By Cauchy--Schwarz, a change of variables in the second factor, and Ingham's fourth moment estimate \cite{Ingham}, the integral is $O_\eps(T^{1+\eps})$. This proves \eqref{eq:eis}.
\end{proof}

\subsection{The oscillatory Bessel transform}\label{subsec:bessel}
To estimate $\mathcal O$, put
\begin{equation}\label{eq:Hxn}
  \mathcal H(x,n)=\mathcal H_{\,g(\cdot)V(n,\cdot)}(x).
\end{equation}
By \eqref{eq:weightbounds} and Lemma~\ref{lem:coeffmean}, the AFE variable in $\mathcal O$ may be restricted to
\[
  n=n_1n_2^2n_3^3n_4^4\ll T^{2+\eps}
\]
with negligible error. We may also restrict the Kloosterman modulus to a fixed power of $T$. Indeed, keeping \eqref{eq:weight} on $\Re w=3$ and shifting the $t$-line in \eqref{eq:besseldef} to
\[
  \Im t=-\left(\frac12-\sigma\right),\qquad 0<\sigma<\frac14,
\]
the power series for the $J$-Bessel function and Stirling's formula give, for some constant $A_\sigma>0$,
\begin{equation}\label{eq:smallx}
  \mathcal H(x,n)\ll_{\sigma,h}T^{A_\sigma}x^{1-2\sigma},
  \qquad 0<x\le1.
\end{equation}
Together with Weil's bound and Lemma~\ref{lem:coeffmean}, this makes the contribution of $c>T^A$ negligible when $A$ is sufficiently large.

\begin{lemma}\label{lem:bessel}
Uniformly for $n\ll T^{2+\eps}$, the transform $\mathcal H(x,n)$ is negligible unless
\begin{equation}\label{eq:besselrange}
  x\gg TM T^{-\eps}.
\end{equation}
In this range, for every fixed $B>0$,
\begin{equation}\label{eq:besselmodel}
  \mathcal H(x,n)=\frac{TM}{\sqrt x}\sum_{\pm}W_\pm(x,n)
       \exp\!\left(\pm i\left(x-\frac{2T^2}{x}\right)\right)
       +O_B(T^{-B}),
\end{equation}
where for fixed $j,k,A\ge0$,
\begin{equation}\label{eq:besselweights}
  (x\partial_x)^j(n\partial_n)^kW_\pm(x,n)
  \ll_{j,k,A,\eps}T^\eps\left(1+\frac n{T^2}\right)^{-A}.
\end{equation}
Consequently, in $\mathcal O$ the modulus is restricted to
\begin{equation}\label{eq:modrange}
  c\ll T^\eps\frac{n_1^{3/2}n_2n_3^{3/2}}{TM}.
\end{equation}
\end{lemma}
\begin{proof}
Using the classical integral representation of the $J$-Bessel function \cite[8.411.11]{GR} and Fourier inversion in $t=T+My$, $\mathcal H(x,n)$ is a linear combination of
\begin{equation}\label{eq:fourierbessel}
  TM\int_\R K_\pm(Mv,n)\exp\{\pm ix\cosh v-2iTv\}\,\dd v,
\end{equation}
where for fixed $j,k,A\ge0$,
\[
  \partial_y^j(n\partial_n)^kK_\pm(y,n)
  \ll_{j,k,A,\eps}T^\eps
  \left(1+\frac n{T^2}\right)^{-A}(1+|y|)^{-A}.
\]
Thus $|Mv|\ll T^\eps$ with negligible error. If $x\ll TM T^{-\eps}$, repeated integration by parts gives \eqref{eq:besselrange}.

For the positive sign, the critical point is
\[
  v_0=\operatorname{arsinh}(2T/x),
\]
with critical value
\[
  \Phi(x,T)=\sqrt{x^2+4T^2}-2T\operatorname{arsinh}(2T/x),
\]
and
\[
  \partial_v^2(x\cosh v-2Tv)|_{v=v_0}=\sqrt{x^2+4T^2}\asymp x.
\]
After the change of variable $y=Mv$, the curvature at the critical point is
\[
  \asymp\frac{x}{M^2}\gg\frac{T}{M}T^{-\eps}\gg T^{\delta-\eps}.
\]
The parameter uniform stationary phase theorem \cite[Main Theorem]{KPY} gives the factor $x^{-1/2}$ in \eqref{eq:besselmodel} and the derivative bounds \eqref{eq:besselweights}. The negative sign is handled in the same way.

Finally,
\[
  \Phi(x,T)=x-\frac{2T^2}{x}+O\!\left(\frac{T^4}{x^3}\right),
  \qquad
  \frac{T^4}{x^3}\ll \frac{T^{1+3\eps}}{M^3}
  \ll T^{-5/4-3\delta+3\eps}.
\]
The exponential of this remainder, together with its fixed Euler derivatives, may be absorbed into $W_\pm$. Taking
\[
  x=\frac{4\pi n_1^{3/2}n_2n_3^{3/2}}c
\]
gives \eqref{eq:modrange}.
\end{proof}

\subsection{Applying the first Poisson summation}\label{subsec:first}
Fix $n_2,n_3,n_4$. Insert a smooth dyadic partition of unity in the $n_1$-sum, and split the $c$-sum into dyadic ranges $c\asymp C$. By the approximate functional equation and Lemma~\ref{lem:bessel}, it is enough to consider
\begin{equation}\label{eq:ranges}
  1\le N_1\le \frac{T^{2+\eps}}{n_2^2n_3^3n_4^4},
  \qquad
  C\ll T^\eps\frac{N_1^{3/2}n_2n_3^{3/2}}{TM}
  \ll \frac{T^{2+\eps}}M.
\end{equation}
For one exponential branch in Lemma~\ref{lem:bessel}, the corresponding block is
\begin{equation}\label{eq:model}
\begin{split}
 \mathscr O(n_2,n_3,n_4;N_1,C)
 ={}&\frac{TM}{n_4^2n_2^{3/2}n_3^{9/4}}
 \sum_{c\asymp C}\frac1{\sqrt c}
 \sum_{n_1\ge1}\frac{S(1,n_2^2n_3^3n_1^3;c)}{n_1^{5/4}}W(n_1,c)\\
 &\qquad\times\e\!\left(
 \frac{2n_2n_3^{3/2}n_1^{3/2}}c
 -\frac{cT^2}{4\pi^2n_2n_3^{3/2}n_1^{3/2}}\right),
\end{split}
\end{equation}
where $W$ is smooth, supported where $N_1\le n_1\le2N_1$, and for all fixed $j,k\ge0$,
\[
  (n_1\partial_{n_1})^j(c\partial_c)^kW(n_1,c)\ll_{j,k,\eps}T^\eps.
\]
The factor $\mu^2(n_3)$ is kept outside the block. The other exponential branch is treated in the same way.

Applying Poisson summation to the innermost $n_1$-sum in \eqref{eq:model}, we obtain
\begin{equation}\label{eq:firstpoisson}
  \mathscr O(n_2,n_3,n_4;N_1,C)
  =\frac{TM}{n_4^2n_2^{3/2}n_3^{9/4}}
   \sum_{c\asymp C}\frac1{c^{3/2}}
   \sum_{m\in\Z}\mathfrak C_c(m;n_2,n_3)\mathcal I_c(m;n_2,n_3),
\end{equation}
where
\begin{equation}\label{eq:Cdef}
  \mathfrak C_c(m;n_2,n_3)
  =\sum_{x\bmod c}S(1,n_2^2n_3^3x^3;c)\e_c(mx),
\end{equation}
and
\begin{equation}\label{eq:firstintegral}
  \mathcal I_c(m;n_2,n_3)=\int_0^\infty x^{-5/4}W(x,c)
  \e\!\left(\frac{2n_2n_3^{3/2}x^{3/2}}c
  -\frac{cT^2}{4\pi^2n_2n_3^{3/2}x^{3/2}}-\frac{mx}c\right)\,\dd x.
\end{equation}

\begin{lemma}\label{lem:first}
The integral $\mathcal I_c(m;n_2,n_3)$ is negligible unless
\[
  m>0,\qquad m\asymp n_2n_3^{3/2}N_1^{1/2}.
\]
In this range,
\begin{equation}\label{eq:firstasymp}
\begin{split}
  \mathcal I_c(m;n_2,n_3)
  ={}&\frac{c^{1/2}}{n_2^{1/2}n_3^{3/4}N_1}
  W_1(m,c)\\
  &\times\e\!\left(-\frac{m^3}{27n_2^2n_3^3c}
  -\frac{27cn_2^2n_3^3T^2}{4\pi^2m^3}\right)
  +O_A(T^{-A}),
\end{split}
\end{equation}
where $W_1$ is smooth, supported in a fixed enlargement of the above range, and for fixed $j,k\ge0$,
\begin{equation}\label{eq:firstweight}
  (m\partial_m)^j(c\partial_c)^kW_1(m,c)\ll_{j,k,\eps}T^\eps.
\end{equation}
\end{lemma}
\begin{proof}
Write the phase in \eqref{eq:firstintegral} as
\[
  \phi_m(x)=\frac{2n_2n_3^{3/2}x^{3/2}}c
  -\frac{cT^2}{4\pi^2n_2n_3^{3/2}x^{3/2}}-\frac{mx}c.
\]
On the support $x\asymp N_1$,
\begin{equation}\label{eq:firstderivative}
  c\phi_m'(x)
  =3n_2n_3^{3/2}\sqrt x
   +\frac{3c^2T^2}{8\pi^2n_2n_3^{3/2}x^{5/2}}-m,
\end{equation}
and the second term on the right is smaller than the first by
\[
  \ll \frac{c^2T^2}{n_2^2n_3^3N_1^3}
  \ll M^{-2}T^{2\eps}
\]
by \eqref{eq:ranges}. Repeated integration by parts therefore gives a negligible contribution unless
$m>0$ and $m\asymp n_2n_3^{3/2}N_1^{1/2}$.

In this range \eqref{eq:firstderivative} has a unique critical point $x_m\asymp N_1$, since
\[
  c\phi_m''(x)
  =\frac{3n_2n_3^{3/2}}{2\sqrt x}
   -\frac{15c^2T^2}{16\pi^2n_2n_3^{3/2}x^{7/2}}
  =\frac{3n_2n_3^{3/2}}{2\sqrt x}
   \{1+O(M^{-2}T^{2\eps})\}.
\]
Expanding the critical point equation gives
\begin{equation}\label{eq:point}
  x_m=\frac{m^2}{9n_2^2n_3^3}
   -\frac{81c^2T^2n_2^2n_3^3}{4\pi^2m^4}
   +O\!\left(\frac{c^4T^4n_2^6n_3^9}{m^{10}}\right).
\end{equation}
Substitution into the phase gives
\begin{equation}\label{eq:dualphase}
  \phi_m(x_m)
  =-\frac{m^3}{27n_2^2n_3^3c}
   -\frac{27cT^2n_2^2n_3^3}{4\pi^2m^3}
   +O\!\left(\frac{c^3T^4n_2^6n_3^9}{m^9}\right).
\end{equation}
The error term in \eqref{eq:dualphase}, together with its fixed Euler derivatives, is
\[
  O\!\left(\frac{T^{1+3\eps}}{M^3}\right),
\]
and its exponential may be absorbed into the weight.

Finally, $\phi_m''(x_m)\asymp n_2n_3^{3/2}/(cN_1^{1/2})$. Since the amplitude in \eqref{eq:firstintegral} has size $N_1^{-5/4}$, stationary phase gives
\[
  N_1^{-5/4}
  \left(\frac{cN_1^{1/2}}{n_2n_3^{3/2}}\right)^{1/2}
  =\frac{c^{1/2}}{n_2^{1/2}n_3^{3/4}N_1}.
\]
The parameter uniform stationary phase theorem \cite[Main Theorem]{KPY} gives \eqref{eq:firstasymp} and the derivative bounds \eqref{eq:firstweight}.
\end{proof}

Applying Lemma~\ref{lem:first} to \eqref{eq:firstpoisson}, we obtain, up to a negligible error,
\begin{equation}\label{eq:postfirst}
\begin{split}
  \mathscr O(n_2,n_3,n_4;N_1,C)
  ={}&\frac{TM}{n_4^2n_2^2n_3^3N_1}
  \sum_{c\asymp C}\frac1c\sum_{m\in\Z}
  \mathfrak C_c(m;n_2,n_3)W_1(m,c)\\
  &\times\e\!\left(-\frac{m^3}{27n_2^2n_3^3c}
  -\frac{27cn_2^2n_3^3T^2}{4\pi^2m^3}\right).
\end{split}
\end{equation}
The smooth weight $W_1$ localizes the $m$-sum to the stationary range stated in Lemma~\ref{lem:first}. Summing \eqref{eq:postfirst} over the smooth dyadic blocks and over $n_2,n_3,n_4$, with the factor $\mu^2(n_3)$, recovers the contribution of the chosen Bessel branch to $\mathcal O$; the other branch is treated identically.

\subsection{Applying the second Poisson summation}\label{subsec:second}
The two oscillatory factors in \eqref{eq:postfirst} now play different roles. The character sum $\mathfrak C_c(m;n_2,n_3)$ is already periodic modulo $c$, while
\[
  \e\!\left(-\frac{m^3}{27n_2^2n_3^3c}\right)
\]
is periodic modulo $27n_2^2n_3^3c$. Thus adjoining the cubic factor enlarges the arithmetic modulus only by the outer factor $27n_2^2n_3^3$. On the effective support $m\asymp n_2n_3^{3/2}N_1^{1/2}$, its phase has size
\[
  \frac{m^3}{27n_2^2n_3^3c}
  \asymp \frac{n_2n_3^{3/2}N_1^{3/2}}{c},
\]
which is of order $TM$ when $c$ is near the upper end of the range in \eqref{eq:ranges}. Hence this large analytic oscillation has been transferred to the arithmetic side by enlarging the modulus only by the outer factor specified above.

By contrast,
\[
  \e\!\left(-\frac{27cn_2^2n_3^3T^2}{4\pi^2m^3}\right)
\]
remains part of the analytic weight. Its phase has size
\[
  \asymp \frac{cT^2}{n_2n_3^{3/2}N_1^{3/2}}
  \ll T^\eps\frac{T}{M},
\]
and is of order $T/M$ near the upper end of the admissible range for $c$. The smaller oscillation in the analytic weight makes the second Poisson summation in $m$ effective.

Applying Poisson summation to the $m$-sum in \eqref{eq:postfirst}, we obtain
\begin{equation}\label{eq:second}
  \frac1{27n_2^2n_3^3c}
  \sum_{\ell\in\Z}\mathfrak S(\ell;n_2,n_3,c)\mathcal J_c(\ell),
\end{equation}
where
\begin{equation}\label{eq:Sdef}
  \mathfrak S(\ell;n_2,n_3,c)
  =\sum_{m\bmod 27n_2^2n_3^3c}\mathfrak C_c(m;n_2,n_3)
    \e_{27n_2^2n_3^3c}(-m^3+\ell m),
\end{equation}
and
\begin{equation}\label{eq:Jdef}
  \mathcal J_c(\ell)=\int_0^\infty W_1(x,c)
   \e\!\left(-\frac{27cT^2n_2^2n_3^3}{4\pi^2x^3}
              -\frac{\ell x}{27n_2^2n_3^3c}\right)\,\dd x.
\end{equation}
Substituting \eqref{eq:second} into \eqref{eq:postfirst}, the dyadic block becomes
\begin{equation}\label{eq:postsecond}
  \mathscr O(n_2,n_3,n_4;N_1,C)
  =\frac{TM}{27n_4^2n_2^4n_3^6N_1}
   \sum_{c\asymp C}\frac1{c^2}
   \sum_{\ell\in\Z}\mathfrak S(\ell;n_2,n_3,c)\mathcal J_c(\ell)
   +O_A(T^{-A}).
\end{equation}
Summing over the smooth dyadic blocks and over $n_2,n_3,n_4$, with the factor $\mu^2(n_3)$, the contribution of the displayed Bessel branch to $\mathcal O$ is therefore, up to a negligible error,
\[
  \frac{TM}{27}
  \sum_{n_2,n_3,n_4\ge1}
  \frac{\mu^2(n_3)}{n_4^2n_2^4n_3^6}
  \sum_{N_1,C}\frac1{N_1}
  \sum_{c\asymp C}\frac1{c^2}
  \sum_{\ell\in\Z}
  \mathfrak S(\ell;n_2,n_3,c)\mathcal J_c(\ell),
\]
where $N_1,C$ satisfy \eqref{eq:ranges}. The other Bessel branch has the same form and is treated identically.

\subsection{The second integral transform}\label{subsec:second-transform}
We now determine the effective range and size of $\mathcal J_c(\ell)$.

\begin{lemma}\label{lem:second}
Let $c\asymp C$. If
\[
  \frac{CT^2}{n_2n_3^{3/2}N_1^{3/2}}\le T^\eps,
\]
then, for every $A>0$,
\[
  \mathcal J_c(\ell)
  \ll_{A,\eps}T^{-A}(1+|\ell|)^{-2}
\]
unless
\[
  |\ell|\ll T^{2\eps}
  \left(1+\frac{n_2n_3^{3/2}C}{N_1^{1/2}}\right),
\]
in which case
\begin{equation}\label{eq:Jsmallpoint}
  |\mathcal J_c(\ell)|
  \ll_\eps T^\eps n_2n_3^{3/2}N_1^{1/2}.
\end{equation}
If
\[
  \frac{CT^2}{n_2n_3^{3/2}N_1^{3/2}}\ge T^\eps,
\]
then, for every $A>0$,
\[
  \mathcal J_c(\ell)
  \ll_{A,\eps}T^{-A}(1+|\ell|)^{-2}
\]
unless
\[
  \ell>0,
  \qquad
  \ell\asymp \frac{C^2T^2}{N_1^2},
\]
in which case
\begin{equation}\label{eq:Jstationary}
  |\mathcal J_c(\ell)|
  \ll_\eps T^\eps
  \frac{n_2^{3/2}n_3^{9/4}N_1^{5/4}}{T\sqrt C}.
\end{equation}
\end{lemma}
\begin{proof}
Put $X=n_2n_3^{3/2}N_1^{1/2}$ and set $x=Xy$ in \eqref{eq:Jdef}. On the support of the resulting weight, $y\asymp1$. Up to absolute constants, the two coefficients in the scaled phase have sizes
\[
  \lambda=\frac{CT^2}{n_2n_3^{3/2}N_1^{3/2}},
  \qquad
  \eta=\frac{|\ell|N_1^{1/2}}{n_2n_3^{3/2}C}.
\]
The derivatives of the weight are $O(T^{\eps_0})$ for every fixed $\eps_0>0$. All parameters occurring here are $T^{O(1)}$ by \eqref{eq:ranges}.

Suppose first that $\lambda\le T^\eps$. If
\[
  |\ell|\gg T^{2\eps}
  \left(1+\frac{n_2n_3^{3/2}C}{N_1^{1/2}}\right),
\]
then $\eta\gg T^{2\eps}$, and the derivative of the scaled phase has size $\gg\eta$ throughout the support. Repeated integration by parts gives, for every fixed $B>0$,
\[
  \mathcal J_c(\ell)
  \ll_{B,\eps_0}XT^{B\eps_0}(1+\eta)^{-B}.
\]
Taking $B$ sufficiently large and then $\eps_0$ sufficiently small, in terms of $A$ and $\eps$, gives
\[
  \mathcal J_c(\ell)
  \ll_{A,\eps}T^{-A}(1+|\ell|)^{-2}.
\]
In the complementary range, the trivial estimate gives
\[
  |\mathcal J_c(\ell)|\ll_\eps T^\eps X,
\]
which is \eqref{eq:Jsmallpoint}.

Now suppose that $\lambda\ge T^\eps$. If $\ell\le0$, or if $\ell>0$ but $\eta$ is outside a fixed enlargement of the range $\eta\asymp\lambda$, the derivative of the scaled phase is bounded below by a fixed multiple of $\lambda+\eta$. Repeated integration by parts, as above, gives
\[
  \mathcal J_c(\ell)
  \ll_{A,\eps}T^{-A}(1+|\ell|)^{-2}.
\]
The remaining range is characterized by
\[
  \ell>0,
  \qquad
  \eta\asymp\lambda,
\]
which is equivalent to
\[
  \ell\asymp \frac{C^2T^2}{N_1^2}.
\]
There the second derivative of the scaled phase is $\asymp\lambda$. The second derivative estimate therefore gives
\[
  |\mathcal J_c(\ell)|
  \ll_\eps T^\eps X\lambda^{-1/2}
  =T^\eps
  \frac{n_2^{3/2}n_3^{9/4}N_1^{5/4}}{T\sqrt C},
\]
which proves \eqref{eq:Jstationary}.
\end{proof}

\subsection{Complete exponential sums}\label{subsec:complete}

We now turn to the arithmetic part of the second Poisson summation.

\begin{proposition}\label{prop:complete-general}
For all $n_2,n_3,c\ge1$ and $\ell\in\Z$,
\begin{equation}\label{eq:complete-general}
  |\mathfrak S(\ell;n_2,n_3,c)|
     \ll_\eps n_2^2n_3^3c^{3/2+\eps}\rho(c)^{1/2},
\end{equation}
where
\[
  \rho(c)=2^{v_2(c)}3^{v_3(c)}
     \prod_{\substack{p\ge5\\v_p(c)\ge2}}p^{v_p(c)}.
\]
Moreover,
\begin{equation}\label{eq:powerful-average}
  \sum_{C<c\le2C}\rho(c)^{1/2}\ll C\log(2C).
\end{equation}
\end{proposition}
\begin{proof}
We begin with the Chinese remainder decomposition. For each prime
$p\mid27n_2^2n_3^3c$, write
\[
  \beta=v_p(n_2^2n_3^3),\qquad
  \nu=v_p(c),\qquad
  j=v_p(27),\qquad
  r=\beta+\nu+j,
\]
and set
\[
  t_p=\left(\frac{27n_2^2n_3^3c}{p^r}\right)^{-1}\pmod{p^r}.
\]
If $\nu\ge1$, also set
\[
  \gamma_p=\left(\frac{c}{p^\nu}\right)^{-1}\pmod{p^\nu}
\]
and define
\begin{equation}\label{eq:Ctwist}
  \mathfrak C_{p^\nu}^{(\gamma_p)}(m;n_2,n_3)
  =\sideset{}{^{*}}\sum_{d\bmod p^\nu}\sum_{x\bmod p^\nu}
  \e_{p^\nu}\!\left(
    \gamma_p\bigl(d+n_2^2n_3^3\bar d x^3+mx\bigr)
  \right).
\end{equation}
When $\nu=0$, we interpret this factor as $1$. The Chinese remainder theorem gives
\begin{equation}\label{eq:S-CRT}
  \mathfrak S(\ell;n_2,n_3,c)
  =\prod_{p\mid27n_2^2n_3^3c}\mathfrak S_p,
\end{equation}
where
\begin{equation}\label{eq:local-S}
  \mathfrak S_p
  =\sum_{m\bmod p^r}
  \mathfrak C_{p^\nu}^{(\gamma_p)}(m;n_2,n_3)
  \e_{p^r}\!\left(t_p(-m^3+\ell m)\right).
\end{equation}
We estimate these local factors case by case.

\medskip
\noindent
\emph{Case 1: $\nu=0$.}
Then
\[
  \mathfrak S_p
  =\sum_{m\bmod p^{\beta+j}}
  \e_{p^{\beta+j}}\!\left(t_p(-m^3+\ell m)\right),
\]
and hence
\begin{equation}\label{eq:local-nu0}
  |\mathfrak S_p|\le p^{\beta+j}\ll p^\beta.
\end{equation}
Here and below the fixed powers of $2$ and $3$ are absorbed into the implied constant.

\medskip
\noindent
\emph{Case 2: $\beta\ge1$ and $\nu\ge1$.}
Thus $p\mid n_2n_3$. If $\nu=1$, then
$S(\gamma_p,\gamma_p n_2^2n_3^3x^3;p)=-1$, and Fourier orthogonality gives
\begin{equation}\label{eq:ramified-C}
  \mathfrak C_p^{(\gamma_p)}(m;n_2,n_3)
  =-p\,\mathbf1_{p\mid m}.
\end{equation}
Therefore
\[
  \mathfrak S_p
  =-p\sum_{\substack{m\bmod p^{\beta+1+j}\\p\mid m}}
  \e_{p^{\beta+1+j}}\!\left(t_p(-m^3+\ell m)\right),
\]
and hence
\begin{equation}\label{eq:ramified-local}
  |\mathfrak S_p|\le p^{\beta+1+j}.
\end{equation}

If $\nu\ge2$, then
\begin{equation}\label{eq:ramified-vanish}
  \mathfrak C_{p^\nu}^{(\gamma_p)}(m;n_2,n_3)=0.
\end{equation}
Indeed, write
\[
  d=d_0+p^{\nu-1}v,
  \qquad v\bmod p,
\]
where $d_0\bmod p^{\nu-1}$ and $(d_0,p)=1$. Then
\[
  \bar d\equiv \bar d_0-p^{\nu-1}v\bar d_0^{\,2}\pmod{p^\nu},
\]
where $\bar d_0$ denotes the inverse modulo $p^\nu$ of the chosen representative $d_0$. Since $p\mid n_2n_3$,
\[
  n_2^2n_3^3\bar d x^3
  \equiv
  n_2^2n_3^3\bar d_0x^3
  \pmod{p^\nu}.
\]
Hence
\[
\begin{split}
  \mathfrak C_{p^\nu}^{(\gamma_p)}(m;n_2,n_3)
  ={}&
  \sideset{}{^{*}}\sum_{d_0\bmod p^{\nu-1}}
  \sum_{x\bmod p^\nu}
  \e_{p^\nu}\!\left(
    \gamma_p\bigl(d_0+n_2^2n_3^3\bar d_0x^3+mx\bigr)
  \right) 
  \sum_{v\bmod p}\e_p(\gamma_pv).
\end{split}
\]
Since $\gamma_p$ is a unit modulo $p$, the last sum vanishes. Thus
\[
  \mathfrak C_{p^\nu}^{(\gamma_p)}(m;n_2,n_3)=0
  \qquad (\nu\ge2),
\]
and consequently $\mathfrak S_p=0$.

\medskip
\noindent
\emph{Case 3: $\beta=0$, $p\ge5$, and $\nu=1$.}
Write $c=pc'$. Then
\[
  \gamma_p=\overline{c'}\pmod p,
  \qquad
  t_p=\overline{27n_2^2n_3^3c'}\pmod p.
\]
In \eqref{eq:Ctwist}, make the changes of variables
$d'=\gamma_pd$ and $x'=\gamma_px$. Dropping the primes afterwards, \eqref{eq:local-S} becomes
\begin{equation}\label{eq:good-local}
\begin{split}
  \mathfrak S_p
  ={}&\sum_{m\bmod p}\sum_{x\bmod p}
  S(1,n_2^2n_3^3c'x^3;p)\e_p(mx) 
  \e_p\!\left(
    \overline{27n_2^2n_3^3c'}(-m^3+\ell m)
  \right).
\end{split}
\end{equation}
Let $\chi_p$ denote the quadratic character modulo $p$. For $x\ne0\bmod p$, the change of variable $d=xd'$ gives
\[
  S(1,n_2^2n_3^3c'x^3;p)
  =S(x,n_2^2n_3^3c'x^2;p).
\]
At $x=0$, the two sides are $-1$ and $p-1$, respectively. Hence the inner Fourier transform in \eqref{eq:good-local} equals
\[
  \sum_{x\bmod p}S(x,n_2^2n_3^3c'x^2;p)\e_p(mx)-p.
\]
Expanding the Kloosterman sum and completing the square in $x$, the first term contributes
\[
  \tau_p\chi_p(n_2^2n_3^3c')
  \sideset{}{^{*}}\sum_{d\bmod p}\sum_{m\bmod p}\chi_p(d)
  \e_p\!\left(
    -\frac{d(d+m)^2}{4n_2^2n_3^3c'}
    +\frac{-m^3+\ell m}{27n_2^2n_3^3c'}
  \right),
\]
where
\[
  \tau_p=\sum_{x\bmod p}\chi_p(x)\e_p(x).
\]
Set $s=3d+m$. For each fixed $d$, the $m$-sum is then equivalently an $s$-sum, and
\[
  -\frac{d(d+m)^2}{4n_2^2n_3^3c'}
  +\frac{-m^3+\ell m}{27n_2^2n_3^3c'}
  =\frac{-s^3+\ell s}{27n_2^2n_3^3c'}
   +\frac{d(3s^2-4\ell)}{36n_2^2n_3^3c'}.
\]
After interchanging the $d$- and $s$-sums, the $d$-sum is the quadratic Gauss sum
\[
  \sideset{}{^{*}}\sum_{d\bmod p}\chi_p(d)
  \e_p\!\left(
    \frac{d(3s^2-4\ell)}{36n_2^2n_3^3c'}
  \right).
\]
Evaluating this sum and using $\tau_p^2=p\chi_p(-1)$ gives
\begin{equation}\label{eq:good-prime-formula}
  \mathfrak S_p
  =p\sum_{s\bmod p}
  \bigl(\chi_p(4\ell-3s^2)-1\bigr)
  \e_p\!\left(
    \overline{27n_2^2n_3^3c'}(-s^3+\ell s)
  \right).
\end{equation}
For the mixed sum in \eqref{eq:good-prime-formula}, Katz \cite[Theorem~1.1(1)]{KatzMixed}, applied in dimension one to
\[
  f(s)=\overline{27n_2^2n_3^3c'}(-s^3+\ell s),
  \qquad
  g(s)=4\ell-3s^2,
\]
with the quadratic character, gives $O(\sqrt p)$, uniformly in $\ell$, including $\ell=0$. The remaining pure cubic sum is also $O(\sqrt p)$ by the Weil bound. Therefore
\begin{equation}\label{eq:good-prime-bound}
  |\mathfrak S_p|\ll p^{3/2}.
\end{equation}

\medskip
\noindent
\emph{Case 4: $\beta=0$ and either $\nu\ge2$ or $p=2,3$.}
The case $\nu=0$ has already been treated, so here $\nu\ge1$. Finite Parseval gives
\[
\begin{split}
  \sum_{m\bmod p^\nu}
  \left|\mathfrak C_{p^\nu}^{(\gamma_p)}(m;n_2,n_3)\right|^2
  &=p^\nu\sum_{x\bmod p^\nu}
    \left|S(\gamma_p,\gamma_p n_2^2n_3^3x^3;p^\nu)\right|^2\\
  &\ll (\nu+1)^2p^{3\nu},
\end{split}
\]
where we used Weil's bound. Hence, by Cauchy--Schwarz,
\begin{equation}\label{eq:Cl1}
  \sum_{m\bmod p^\nu}
  \left|\mathfrak C_{p^\nu}^{(\gamma_p)}(m;n_2,n_3)\right|
  \ll (\nu+1)p^{2\nu}.
\end{equation}
Since this function is periodic modulo $p^\nu$, while the $m$-sum in \eqref{eq:local-S} is modulo $p^{\nu+j}$, we obtain
\begin{equation}\label{eq:bad-prime-bound}
  |\mathfrak S_p|\ll p^j(\nu+1)p^{2\nu}.
\end{equation}

We now multiply the local estimates in \eqref{eq:S-CRT}. The primes with $\nu=0$ contribute at most the corresponding factor of $n_2^2n_3^3$. At primes $p\ge5$ with $\nu=1$, \eqref{eq:good-prime-bound} gives the square root bound $p^{3/2}$. In Case~4,
\[
  p^j(\nu+1)p^{2\nu}
  \ll_\eps p^{3\nu/2+\eps\nu}p^{\nu/2},
\]
and the additional factor $p^{\nu/2}$ is precisely recorded by $\rho(c)^{1/2}$. Case~2 is no larger, by \eqref{eq:ramified-local} and \eqref{eq:ramified-vanish}. Absorbing the factors $\nu+1$ into $c^\eps$, we obtain
\[
  |\mathfrak S(\ell;n_2,n_3,c)|
  \ll_\eps n_2^2n_3^3c^{3/2+\eps}\rho(c)^{1/2},
\]
which proves \eqref{eq:complete-general}.

It remains to prove \eqref{eq:powerful-average}. Put $t=\rho(c)$. Since $t\mid c$, for each fixed $t\le2C$ there are $O(C/t)$ integers $c$ with
$C<c\le2C$ and $\rho(c)=t$. Moreover, $t$ has the form
\[
  t=2^a3^b u^2v^3,
\]
where $a,b\ge0$, $v$ is squarefree, and $(uv,6)=1$. Indeed, for every exponent $e\ge2$, write $e=2r$ if $e$ is even and $e=2r+3$ if $e$ is odd; hence the part of $t$ supported on primes $p\ge5$ has the form $u^2v^3$ with $v$ squarefree. Therefore
\[
\begin{split}
  \sum_{C<c\le2C}\rho(c)^{1/2}
  &\ll C\sum_{t\le2C}t^{-1/2}\\
  &\ll C
  \sum_{a,b\ge0}2^{-a/2}3^{-b/2}
  \sum_{u\le\sqrt{2C}}\frac1u
  \sum_{v\ge1}\frac1{v^{3/2}}\\
  &\ll C\log(2C),
\end{split}
\]
as required.
\end{proof}

\subsection{Estimate for the Kloosterman term}\label{subsec:summing}
\begin{proposition}\label{prop:block}
For the dyadic block in \eqref{eq:model},
\begin{equation}\label{eq:blockbound}
  \mathscr O(n_2,n_3,n_4;N_1,C)
  \ll_{h,\delta,\eps}T^\eps
  \left(
    \frac{MN_1^{1/4}}{n_4^2n_2^{1/2}n_3^{3/4}}
    +
    \frac{n_2^{3/2}n_3^{9/4}N_1^{5/4}}{n_4^2M}
  \right).
\end{equation}
\end{proposition}
\begin{proof}
By Proposition~\ref{prop:complete-general} and \eqref{eq:powerful-average}, \eqref{eq:postsecond} gives
\begin{equation}\label{eq:master}
  \mathscr O(n_2,n_3,n_4;N_1,C)
  \ll_\eps
  T^\eps\frac{TM\sqrt C}{n_4^2n_2^2n_3^3N_1}
  \sup_{c\asymp C}\sum_{\ell\in\Z}|\mathcal J_c(\ell)|.
\end{equation}

Suppose first that
\[
  \frac{CT^2}{n_2n_3^{3/2}N_1^{3/2}}\le T^\eps.
\]
Summing the bounds in Lemma~\ref{lem:second} over $\ell$, and using \eqref{eq:ranges}, gives
\[
  \sum_{\ell\in\Z}|\mathcal J_c(\ell)|
  \ll_\eps T^\eps n_2n_3^{3/2}N_1^{1/2}.
\]
Moreover,
\[
  \sqrt C\ll T^\eps
  \frac{n_2^{1/2}n_3^{3/4}N_1^{3/4}}T.
\]
Hence \eqref{eq:master} gives
\[
  \mathscr O(n_2,n_3,n_4;N_1,C)
  \ll_\eps T^\eps
  \frac{MN_1^{1/4}}{n_4^2n_2^{1/2}n_3^{3/4}}.
\]

Suppose now that
\[
  \frac{CT^2}{n_2n_3^{3/2}N_1^{3/2}}\ge T^\eps.
\]
Lemma~\ref{lem:second} gives
\[
  \sum_{\ell\in\Z}|\mathcal J_c(\ell)|
  \ll_\eps T^\eps
  \frac{n_2^{3/2}n_3^{9/4}N_1^{5/4}}{T\sqrt C}
  \left(1+\frac{C^2T^2}{N_1^2}\right).
\]
Substituting this into \eqref{eq:master}, and using
\[
  \frac{C^2T^2}{N_1^2}
  \ll T^\eps\frac{n_2^2n_3^3N_1}{M^2}
\]
from \eqref{eq:ranges}, gives
\[
  \mathscr O(n_2,n_3,n_4;N_1,C)
  \ll_\eps T^\eps
  \left(
    \frac{MN_1^{1/4}}{n_4^2n_2^{1/2}n_3^{3/4}}
    +
    \frac{n_2^{3/2}n_3^{9/4}N_1^{5/4}}{n_4^2M}
  \right).
\]
This proves \eqref{eq:blockbound}.
\end{proof}

\begin{proposition}\label{prop:offbound}
The Kloosterman term satisfies
\begin{equation}\label{eq:offbound}
  \mathcal O\ll_{h,\delta,\eps} \frac{T^{5/2+\eps}}M.
\end{equation}
\end{proposition}
\begin{proof}
By \eqref{eq:ranges},
\[
  N_1^{1/2}\ll \frac{T^{1+\eps}}
  {n_4^2n_2n_3^{3/2}}.
\]
Hence the two terms in Proposition~\ref{prop:block} are respectively
\[
  \ll
  \frac{MT^{1/2+\eps}}{n_4^3n_2n_3^{3/2}},
  \qquad
  \frac{T^{5/2+\eps}/M}{n_4^7n_2n_3^{3/2}}.
\]
The sums over $n_3$ and $n_4$ converge absolutely, while the $n_2$-sum costs only a logarithm. The dyadic sums in $N_1$ and $C$ also cost logarithms. Absorbing these into $T^\eps$, and including both exponential branches, proves \eqref{eq:offbound}.
\end{proof}

\subsection{Completion of the proof}\label{subsec:completion}

\begin{proof}[Proof of Theorem~\ref{thm:main}]
Combining \eqref{eq:decomposition}, Lemma~\ref{lem:diag-eis}, and Proposition~\ref{prop:offbound}, we obtain
\[
  \mathcal M_g
  =\frac13 TM\int_\R h(x)\,\dd x
  +O_{h,\delta,\eps}\!\left( \frac{T^{5/2+\eps}}M
    \right),
\]
since $T^{1+\eps}+MT^\eps\ll MT^{1/2+\eps}\ll T^{5/2+\eps}/M$ in the stated range.

For $t_\phi\ge0$, the second term in
\[
  g(t_\phi)
  =h\!\left(\frac{t_\phi-T}{M}\right)
   +h\!\left(\frac{t_\phi+T}{M}\right)
\]
is negligible by the rapid decay of $h$, since $T/M\ge T^\delta$, together with the weighted Weyl law. Hence the same asymptotic formula holds with $g(t_\phi)$ replaced by $h((t_\phi-T)/M)$, proving \eqref{eq:main}.
\end{proof}

\section{Applications of Theorem~\ref{thm:main}}\label{sec:applications}

\subsection{Quantitative non-vanishing}
\begin{proof}[Proof of Corollary~\ref{cor:nonvan}]
To localize the lower bound to the actual interval rather than a Gaussian tail, take
\[
  h_*(v)=e^{-4v^2}-e^{-3}e^{-v^2}.
\]
It satisfies
\[
  \int_{\R} h_*(v)\,\dd v=\sqrt\pi(1/2-e^{-3})>0,
\]
is nonpositive for $|v|\ge1$, and is at most $1$. Theorem~\ref{thm:main} and positivity give
\[
  \sum_{|t_\phi-T|\le M}
    \frac{L(1/2,\sym^3\phi)}{L(1,\sym^2\phi)}\gg TM.
\]
The bounds of Hoffstein--Lockhart \cite{HLock} and Molteni \cite{Molteni} imply
$T^{-\eps}\ll L(1,\sym^2\phi)^{-1}\ll T^\eps$ on the window.
It follows that
\begin{equation}\label{eq:firstunweighted}
  \sum_{\substack{|t_\phi-T|\le M\\L(1/2,\sym^3\phi)\ne0}}
  L(1/2,\sym^3\phi)\gg_\eps TM T^{-\eps}.
\end{equation}
Huang--Li \cite[Theorem~1.4]{HL} prove
\begin{equation}\label{eq:tenth}
  \sum_{T\le t_\phi\le2T}|L(1/2,\sym^3\phi)|^{10}\ll_\eps T^{10+\eps}.
\end{equation}
The same bound holds for the sum over $|t_\phi-T|\le M$: for sufficiently large $T$, we have $M\le T/2$, so we may apply \eqref{eq:tenth} on $[T/2,T]$ and $[T,2T]$.
Applying H\"older's inequality with exponents $10/9$ and $10$, we obtain
\[
  \sum_{\substack{|t_\phi-T|\le M\\L(1/2,\sym^3\phi)\ne0}}
  L(1/2,\sym^3\phi)
  \le
  \bigg(
    \sum_{\substack{|t_\phi-T|\le M\\L(1/2,\sym^3\phi)\ne0}}1
  \bigg)^{9/10}
  \bigg(
    \sum_{|t_\phi-T|\le M}|L(1/2,\sym^3\phi)|^{10}
  \bigg)^{1/10}.
\]
Combining this with \eqref{eq:firstunweighted} and \eqref{eq:tenth} gives
\[
  \#\{\phi:|t_\phi-T|\le M,\ L(1/2,\sym^3\phi)\ne0\}
  \gg_{\delta,\eps} M^{10/9}T^{-\eps},
\]
as claimed.
\end{proof}

\subsection{Almost-all Lindel\"of bounds}
\begin{proof}[Proof of Corollary~\ref{cor:lindelof}]
Use a positive Gaussian in Theorem~\ref{thm:main}, bounded below on $|t-T|\le M$. Positivity from Lemma~\ref{lem:fe} gives
\[
  \sum_{\substack{|t_\phi-T|\le M\\L(1/2,\sym^3\phi)>Y}}
      \frac1{L(1,\sym^2\phi)}\ll\frac{TM}{Y}.
\]
Since $L(1,\sym^2\phi)^{-1}\gg_\eps T^{-\eps}$ on the window, we obtain
\[
  \#\{\phi:|t_\phi-T|\le M,\ L(1/2,\sym^3\phi)>Y\}
  \ll_{\delta,\eps}\frac{TM T^\eps}{Y}.
\]
Taking $Y=T^a$ and choosing $\eps>0$ sufficiently small gives
\[
  \#\{\phi:|t_\phi-T|\le M,\ L(1/2,\sym^3\phi)>T^a\}
  \ll_{\delta,a} TM T^{-a/2},
\]
which proves the claimed almost-all Lindel\"of statement.
\end{proof}

\subsection{Almost-all cubic periods}
\enlargethispage{2\baselineskip}
\begin{proof}[Proof of Corollary~\ref{cor:cubic}]
Recall from Lemma~\ref{lem:fe} that $L(1/2,\sym^3\phi)\ge0$.
The ordinary central value is also nonnegative: for even $\phi$ this is due to Katok--Sarnak \cite[Corollary~1]{KatokSarnak}, while for odd $\phi$ the functional equation has root number $-1$, so $L(1/2,\phi)=0$.

Watson's triple product formula \cite{Watson} 
and Stirling's formula give, for $t_\phi\asymp T$,
\[
  |\langle1,\phi^3\rangle|
  \ll_\eps T^{-1+\eps}
  \frac{L(1/2,\sym^3\phi)^{1/2}
        L(1/2,\phi)}
       {L(1,\sym^2\phi)^{3/2}}.
\]
Using $L(1,\sym^2\phi)^{-1}\ll_\eps T^\eps$ and Cauchy's inequality,
\[
\begin{split}
  \sum_{T\le t_\phi\le T+M}|\langle1,\phi^3\rangle|
  \ll_\eps T^{-1+\eps}
  &\left(\sum_{T\le t_\phi\le T+M}
    \frac{L(1/2,\sym^3\phi)}{L(1,\sym^2\phi)}\right)^{1/2}
    \left(\sum_{T\le t_\phi\le T+M}
    \frac{L(1/2,\phi)^2}{L(1,\sym^2\phi)}\right)^{1/2}.
\end{split}
\]
The first sum is $O_{\delta,\eps}(TM T^\eps)$ by Theorem~\ref{thm:main}, applied with a nonnegative Gaussian bounded below on the interval. For the second sum, the approximate functional equation and the local spectral large sieve inequality (see, e.g., Jutila \cite[Theorem~1.1]{Jutila00}) give
\[
  \sum_{T\le t_\phi\le T+M}
    \frac{L(1/2,\phi)^2}{L(1,\sym^2\phi)}
  \ll_\eps TM T^\eps.
\]
This proves \eqref{eq:cubic-first}. For each fixed $a>0$, Markov's inequality and the choice $\eps=a/2$ give
\[
\begin{split}
  &\#\{\phi:T\le t_\phi\le T+M,\ |\langle1,\phi^3\rangle|>T^{-1+a}\}\\
  &\qquad\le T^{1-a}\sum_{T\le t_\phi\le T+M}|\langle1,\phi^3\rangle|
  \ll_{\delta,a}TM T^{-a/2}.
\end{split}
\]
This proves the almost-all statement.
\end{proof}

\section{The holomorphic case}\label{sec:holo}

\subsection{Approximate functional equation and diagonal}
Recall the notation of Theorem~\ref{thm:holo}. In this section, $M$ denotes the length of the weight window. Put
\begin{equation}\label{eq:dk}
  d_k=\psi(k/2)+\psi((3k-2)/2)-2\log(2\pi),\qquad \psi=\Gamma'/\Gamma.
\end{equation}
The archimedean factor and the functional equation are
\begin{equation}\label{eq:hologamma}
 \begin{gathered}
  \gamma_k(s)=\GC\!\left(s+\frac{k-1}{2}\right)
                 \GC\!\left(s+\frac{3(k-1)}2\right),\qquad
  \GC(s)=2(2\pi)^{-s}\Gamma(s),\\
  \Lambda(s,\sym^3 f):=\gamma_k(s)L(s,\sym^3f)
  =-\Lambda(1-s,\sym^3 f).
 \end{gathered}
\end{equation}
Indeed, the weight $k$ discrete series Weil parameter $\rho_{k-1}$ satisfies
$\sym^3\rho_{k-1}=\rho_{3(k-1)}\oplus\rho_{k-1}$, whose epsilon factor is
$i^{3k-2}i^k=-1$; the finite epsilon factors are $1$. See
\cite[Theorem~4.2]{KS00}. In particular, $L(1/2,\sym^3 f)=0$.

Write
\[
  L(s,\sym^3 f)=\sum_{n\ge1}\frac{A_f(n)}{n^s}.
\]
Consider
\[
  \intc{3}
  \frac{\gamma_k(1/2+w)L(1/2+w,\sym^3 f)}{\gamma_k(1/2)}
  \frac{e^{w^2}}{w^2}\,\dd w.
\]
Shifting the contour and using \eqref{eq:hologamma} in the standard way gives
\begin{equation}\label{eq:holoafe}
 \begin{gathered}
  L'(1/2,\sym^3 f)=2\sum_{n\ge1}\frac{A_f(n)}{\sqrt n}V_k^\dagger(n),\\
  V_k^\dagger(n)=\intc{3}\frac{e^{w^2}}{w^2}
             \frac{\gamma_k(1/2+w)}{\gamma_k(1/2)}n^{-w}\,\dd w.
 \end{gathered}
\end{equation}
Uniformly for $k\asymp K$, and for all fixed $j,b,B\ge0$,
\begin{equation}\label{eq:holo-weight-bounds}
  (n\partial_n)^j\partial_k^b V_k^\dagger(n)
  \ll_{j,b,B,\eps} K^{-b+\eps}(1+n/K^2)^{-B}.
\end{equation}
Moving the contour in \eqref{eq:holoafe} to the left, the residue at $w=0$ is
$d_k-\log n$, and
\begin{equation}\label{eq:dkasym}
  d_k=2\log k+\log\!\left(\frac3{16\pi^2}\right)+O(k^{-1}).
\end{equation}
Thus, the derivative introduces a logarithmic weight, while the natural length remains $k^2$.

The Kuznetsov formula is replaced by the Petersson formula. With the usual Hecke normalization,
\begin{equation}\label{eq:petersson}
 \begin{gathered}
  \omega_f=\frac{2\pi^2}{(k-1)L(1,\sym^2 f)},\\
  \sum_{f\in H_k}\omega_f\lambda_f(m)\lambda_f(n)
   =\delta_{mn}+2\pi i^k\sum_{c\ge1}\frac{S(m,n;c)}c
         J_{k-1}\!\left(\frac{4\pi\sqrt{mn}}c\right).
 \end{gathered}
\end{equation}
The coefficient identity of Lemma~\ref{lem:coeffidentity} applies unchanged. Hence, the diagonal again comes from $n_1=n_2=n_3=1$ and, for fixed $k$, equals
\begin{equation}\label{eq:holodiag}
  2\intc{3}\frac{e^{w^2}}{w^2}
     \frac{\gamma_k(1/2+w)}{\gamma_k(1/2)}\zeta(2+4w)\,\dd w
   =2\zeta(2)d_k+8\zeta'(2)+O(k^{-1/2+\eps}).
\end{equation}
For the second equality, shift to $\Re w=-1/4+\sigma$, with $0<\sigma<1/4$, and take the double pole at $w=0$. Multiplying by $(k-1)/(2\pi^2)$ and averaging over $k$ gives precisely the main sum in \eqref{eq:holomain}, with error $O(MK^{1/2+\eps})$.

\subsection{Kloosterman terms and completion}
It remains to estimate the Kloosterman terms arising from the Petersson formula. The role of the Kuznetsov Bessel transform in Section~\ref{sec:proof} is now played by
\[
  \frac2\pi
  \sum_{\substack{k\ge2\\ k\equiv0\bmod 2}}
  (k-1)\Psi\!\left(\frac{k-K}{M}\right)
  V_k^\dagger(n)i^kJ_{k-1}(x).
\]
Using the integral representation of $J_{k-1}$, Poisson summation in the even weight $k$, and stationary phase, as in \cite[Proposition~8.1 and Corollary~8.2]{ILS}, one obtains the holomorphic analogue of Lemma~\ref{lem:bessel}: the transform is negligible unless
\[
  x\gg KMK^{-\eps},
\]
and in this range it has the form
\begin{equation}\label{eq:holobessel-short}
  \frac{KM}{\sqrt x}\sum_{\pm}W_\pm(x,n)
  \exp\!\left(\pm i\left(x+\frac{(K-1)^2}{2x}\right)\right)
  +O_A(K^{-A}),
\end{equation}
where the weights $W_\pm$ satisfy the same type of derivative bounds as in \eqref{eq:besselweights}, with $T$ replaced by $K$. The bounds in \eqref{eq:holo-weight-bounds} ensure the required smoothness in the $k$-variable; there is no need to remove the dependence of $V_k^\dagger(n)$ on $k$.

From this point the argument is the same as in Section~\ref{sec:proof}. After inserting the coefficient identity, applying the first Poisson summation in $n_1$, and performing stationary phase, the same cubic arithmetic factor
\[
  \e\!\left(-\frac{m^3}{27n_2^2n_3^3c}\right)
\]
appears. The reciprocal analytic phase has the opposite sign from the Maass case, which only reverses the sign of the stationary second dual variable. The second Poisson summation therefore produces exactly the same complete sums $\mathfrak S(\ell;n_2,n_3,c)$ as in \eqref{eq:Sdef}, and Proposition~\ref{prop:complete-general} applies without change. Repeating the estimates of Lemma~\ref{lem:second} and Propositions~\ref{prop:block}--\ref{prop:offbound}, with $T$ replaced by $K$, gives
\begin{equation}\label{eq:holooff}
  \mathcal O_{\mathrm{hol}}(K,M)
  \ll_{\Psi,\delta,\eps}
  \frac{K^{5/2+\eps}}M.
\end{equation}

\begin{proof}[Proof of Theorem~\ref{thm:holo}]
Combining \eqref{eq:holodiag} and \eqref{eq:holooff} proves \eqref{eq:holomain}. By \eqref{eq:dkasym} and the spacing $2$ of the even weights, the leading term is
\[
  \frac16 KM\log K\int_\R\Psi(x)\,\dd x.
\]
If $\Psi$ is supported in $[-1,1]$ and $\int_\R\Psi(x)\,\dd x>0$, then the main term dominates the error term for $K^{3/4+\delta}\le M\le K^{1-\delta}$. Hence the weighted first moment is nonzero for all sufficiently large $K$, and therefore $L'(1/2,\sym^3 f)\ne0$ for at least one even weight $|k-K|\le M$.
\end{proof}

\section*{Acknowledgments}
The authors thank Jianya Liu, Zeev Rudnick, and Peter Sarnak for their interest and encouragement.
This material is based upon work supported by the Swedish Research Council under grant no. 2021-06594 while the authors were in residence at the Institut Mittag-Leffler in Djursholm, Sweden, during the program ``Analytic Number Theory'' (January to April 2024), when the initial work on this project was carried out. We are grateful to Institut Mittag-Leffler for providing excellent working conditions.

\section*{Use of AI tools}
During the development of this work, the authors used ChatGPT to help derive Lemma 2.3, carry out some auxiliary calculations, and improve the presentation of the manuscript. 
The authors take full responsibility for the mathematical content of this paper, including all results and proofs.
%All mathematical results and proofs presented in this paper have been independently verified by the authors, who take full responsibility for their correctness.

\end{document}